\documentclass[12 pt]{article}
\usepackage[centertags]{amsmath}
\usepackage{amssymb}
\usepackage{amsthm}
\usepackage{newlfont}
\usepackage{amsfonts}

\usepackage{graphicx}
\usepackage{color}
\usepackage[colorlinks]{hyperref}
\usepackage{mathptmx}      
\begin{document}

\begin{center}\Large \textbf{Erratum to: "Remarks on the second neighborhood problem''}
\end{center}

\begin{center}
  Salman Ghazal\footnote{Department of Mathematics, Faculty of Sciences I, Lebanese University, Hadath, Beirut, Lebanon.\\
                       E-mail: salmanghazal@hotmail.com\\
                       \par Institute Camille Jordan, D\'{e}partement de Math\'{e}matiques, Universit\'{e} Claude Bernard Lyon 1, France.
                       }
\end{center}

\begin{abstract}

We prove that the proof of existence of weighted local median order of weighted tournaments is wrong and that the proof of the correct statement which asserts that every digraph obtained from a tournament by deleting a set of arcs incident to the same vertex contains a mistake, in the paper entitled "Remarks on the second neighborhood problem". We introduce correct proofs of each.

\end{abstract}

 Let $T=(V, E)$ be a tournament on finite number of vertices. Let $L=v_1v_2...v_n$ be an order of $V$. $L$ is called a median order of $T$ if $\omega (L) :=|\{(v_i,v_j)\in E(D); i<j\}|$ is maximized. The interval $[i,j]$ is the sub-digraph induced by the set $\{v_i,v_{i+1}, ...,v_j\}$. $L$ is called a local median order of $T$ if it satisfies the feedback property: For all $1\leq i\leq j\leq n:$
$$ |N^{+}_{[i,j]}(v_i) | \geq  |N^{-}_{[i,j]}(v_i) | $$
and
$$ |N^{-}_{D[i,j]}(v_j) | \geq  | N^{+}_{D[i,j]}(v_j) | .$$

Every median order is a local median order.\\

 Indeed, let $L=v_1v_2...v_n$ be a median order of $T$ and assume that the feedback property does not hold for some $i\leq j$. Suppose that $ |N^{+}_{[i,j]}(v_i) | <  |N^{-}_{[i,j]}(v_i)|  $. Let $L'=v_1...v_{i-1}v_{i+1}...v_jv_iv_{j+1}...v_n$ be the order obtained from $L$ by inserting $v_i$ just after $v_j$. Then we have:

$$\omega(L')=\omega(L)+|\{(v_k,v_i)\in E(D); i\leq k\leq j\}|-|\{(v_i,v_k)\in E(D); i\leq k\leq j\}|$$ $$=\omega(L)+|N^{-}_{[i,j]}(v_i)|-|N^{+}_{[i,j]}(v_i) | >\omega(L)$$ which contradicts the maximality of $\omega(L)$.\\

Suppose that $|N^{-}_{D[i,j]}(v_j) | <  | N^{+}_{D[i,j]}(v_j) |$. Let $L''=v_1...v_{i-1}v_jv_iv_{i+1}...v_{j-1}v_{j+1}...v_n$ be the order obtained from $L$ by inserting $v_j$ just before $v_i$. Then we have:

$$\omega(L'')=\omega(L)+|\{(v_j,v_k)\in E(D); i\leq k\leq j\}|-|\{(v_k,v_j)\in E(D); i\leq k\leq j\}|$$ $$=\omega(L)+|N^{+}_{[i,j]}(v_j)|-|N^{-}_{[i,j]}(v_j) | >\omega(L)$$ which contradicts the maximality of $\omega(L)$. Therefore, $L$ is a local median order of $T$ (see \cite{m.o.}).\\

Now we consider the case of weighted tournaments, that is there is a positive weight mapping $\omega: V \longrightarrow \mathbb{R}^{+*} $ that assigns to each vertex $x\in V$ a weight $\omega (x)>0$. In \cite{fidler}, an order $L=v_1v_2...v_n$ of $V$ is called a weighted local median order of $T$ if it satisfies the feedback property: For all $1\leq i\leq j\leq n:$
$$ \displaystyle\sum_{x\in N^{+}_{[i,j]}(v_i) }\omega (x) \geq  \sum_{x\in N^{-}_{[i,j]}(v_i)}\omega (x)  $$
and
$$ \displaystyle\sum_{x\in N^{-}_{D[i,j]}(v_j) }\omega (x) \geq  \sum_{x\in N^{+}_{D[i,j]}(v_j)} \omega (x) .$$

 In \cite{fidler}, the weight of an arc $e=(u,v)\in E$ is the weight of its tail, that is, $\omega (e)=\omega (u)$. Let $L=v_1v_2...v_n$ be an order of $V$. The weight of $L$ is defined in \cite{fidler} as $\omega (L):=\displaystyle\sum_{(v_i,v_j)\in E, \, i<j}\omega((v_i,v_j))$.
Suppose that $L$ is an order of $T$ whose weight is maximum. The authors in \cite{fidler}, claimed that $L$ is a local median order of $T$ without proving it. If we proceed to prove their claim as in the proof of the non-weighted case \cite{m.o.}, then the proof will fail. Following the proof of the non weighted case,  we let $L=v_1v_2...v_n$ be an order of $T$ with maximum weight and assume that the feedback property does not hold for some $i\leq j$. We suppose that $ \displaystyle\sum_{x\in N^{+}_{[i,j]}(v_i) }\omega (x) <  \sum_{x\in N^{-}_{[i,j]}(v_i)}\omega (x)  $. Let $L'=v_1...v_{i-1}v_{i+1}...v_jv_iv_{j+1}...v_n$ be the order obtained from $L$ by inserting $v_i$ just after $v_j$. Then we have:

$$\displaystyle\omega (L')= \omega (L)+\sum_{(v_k,v_i)\in E, \, i\leq k\leq j}\omega((v_k,v_i))-\sum_{(v_i,v_k)\in E, \, i\leq k\leq j}\omega((v_i,v_k))$$
$$\displaystyle= \omega (L)+\sum_{(v_k,v_i)\in E, \, i\leq k\leq j}\omega(v_k)-\sum_{(v_i,v_k)\in E, \, i\leq k\leq j}\omega(v_k)$$
$$=\omega (L)+ \displaystyle\sum_{x\in N^{-}_{[i,j]}(v_i) }\omega (x) - \,  |N^{+}_{[i,j]}(v_i)|.\omega(v_i)$$ which is not necessarily equal to
$$\omega(L)+\displaystyle\sum_{x\in N^{-}_{[i,j]}(v_i) }\omega (x)-\sum_{x\in N^{-}_{[i,j]}(v_i)}\omega (x) .$$ So, a contradiction is not reached. \\
 Similarly, if we suppose that $ \displaystyle\sum_{x\in N^{-}_{D[i,j]}(v_j) }\omega (x) <  \sum_{x\in N^{+}_{D[i,j]}(v_j)} \omega (x) $, a contradiction is not reached.  Let $L''=v_1...v_{i-1}v_jv_iv_{i+1}...v_{j-1}v_{j+1}...v_n$ be the order obtained from $L$ by inserting $v_j$ just before $v_i$. Then we have:

$$\displaystyle \omega(L'')=\omega(L)+|N^{+}_{[i,j]}(v_j)|.\omega(v_i)-\, \sum_{x\in N^{-}{[i,j]}(v_i) }\omega (x)$$ which is not necessarily equal to  $$\omega(L)+\displaystyle\sum_{x\in N^{+}_{[i,j]}(v_j) }\omega (x)-\, \sum_{x\in N^{-}_{[i,j]}(v_j)}\omega (x) .$$

However, we prove the existence of weighted local median as follows: We define the weight of an arc $e=(u,v)\in E$ as $\omega (e)=\omega (u). \omega(v)$. Let $L=v_1v_2...v_n$ be an order of $V$. We define the weight of $L$ as $\omega (L):=\displaystyle\sum_{(v_i,v_j)\in E, \, i<j}\omega((v_i,v_j))$. Hence $\omega (L):=\displaystyle\sum_{(v_i,v_j)\in E, \, i<j}\omega(v_i).\omega(v_j)$. An order with maximum weight is called a weighted median order of $T$.\\
 Now, every weighted median order is a weighted local median order. Indeed, let $L=v_1v_2...v_n$ be a weighted median order of $T$ and assume that the feedback property does not hold for some $i\leq j$. Suppose that $ \displaystyle\sum_{x\in N^{+}_{[i,j]}(v_i) }\omega (x) <  \sum_{x\in N^{-}_{[i,j]}(v_i)}\omega (x)  $. Let $L'=v_1...v_{i-1}v_{i+1}...v_jv_iv_{j+1}...v_n$ be the order obtained from $L$ by inserting $v_i$ just after $v_j$. Then we have:

$$\displaystyle\omega (L')= \omega (L)+\sum_{(v_k,v_i)\in E, \, i\leq k\leq j}\omega((v_k,v_i))-\sum_{(v_i,v_k)\in E, \, i\leq k\leq j}\omega((v_i,v_k))$$
$$\displaystyle= \omega (L)+\sum_{(v_k,v_i)\in E, \, i\leq k\leq j}\omega(v_k).\omega(v_i)-\sum_{(v_i,v_k)\in E, \, i\leq k\leq j}\omega(v_k).\omega(v_i)$$
$$=\omega (L)+ \displaystyle[\sum_{x\in N^{-}_{[i,j]}(v_i) }\omega (x) - \sum_{x\in N^{+}_{[i,j]}(v_i) }\omega (x)].\omega(v_i) \, >\omega(L).$$ which is a contradiction to the maximality of $\omega(L)$. \\
 Similarly, suppose that $ \displaystyle\sum_{x\in N^{-}_{D[i,j]}(v_j) }\omega (x) <  \sum_{x\in N^{+}_{D[i,j]}(v_j)} \omega (x) $.  Let $L''=v_1...v_{i-1}v_jv_iv_{i+1}...v_{j-1}v_{j+1}...v_n$ be the order obtained from $L$ by inserting $v_j$ just before $v_i$. Then we have:

$$\displaystyle \omega(L'')=\omega (L)+ \displaystyle[\sum_{x\in N^{+}_{[i,j]}(v_j) }\omega (x) - \sum_{x\in N^{-}_{[i,j]}(v_j) }\omega (x)].\omega(v_j) \, >\omega(L),$$ which is a contradiction. So, $L$ satisfies the feedback property and thus it is a local median order of $T$.

Although the statement of theorem 3.5 of \cite{fidler} is true, but its proof contains a mistake. The theorem states that if $G=(V,E)$ is an orientation obtained from a tournament by deleting a star (set of arcs incident with the same vertex), then $G$ satisfies Seymour's second neighborhood conjecture.
Their proofs proceeds as follows. Let $x$ denote the vertex incident with all the missing edges ($x$ is the center of the deleted star). For every missing edge $xy$, let $Q=\{q\in V; (y,q), (q,x), y\notin N^{++}(q)\}$ and $R=\{r\in V; (r,y)\in E, (x,r)\notin E, x\notin N^{++}(r)\}$. Then at least one of $Q$ or $R$ is empty. Indeed, assume that $q\in Q$ and $r\in R$. Since $q$ and $r$ are adjacent vertices in $G$, then there are two possibilities. If $(r,q)\in E$, then $x\in N^{++}(r)$, a contradiction. If $(q,r)\in E$, then $y\in N^{++}(q)$, a contradiction. Assume first that $Q=\phi$. In this case add the arc $(x,y)$ to $G$. No vertex, except possibly $x$, received a new vertex in its second out-neighborhood. Assume next that $R=\phi$. In this case, add the arc $(y,x)$. \\
Now the authors of \cite{fidler} claims that no vertex at all received a new vertex in its second out-neighborhood and that the obtained tournament $T$, by adding all these new arcs using the above procedure, satisfies the following: $\forall v\in V-x$,  $N^{++}(v)=N^{++}_T(v)$. This is false. As a counterexample, consider the following digraph $G$ with vertex set $\{x,y,z,t,q\}$ and arc set $\{(y,q), (q,x), (z,q), (x,z), (z,y), (t,q),(t,z),(t,y)\}$. Then $q\in Q$, $R=\phi$ and $z\notin N^{++}(y)$. So we must add the arc $(y,x)$ to $G$. But now $z$ becomes a new second out-neighbor of $y$. \\

 A correct proof can be established as follows. We orient each missing edge $xy$ by following the above procedure exactly to obtain a tournament $T$. If $v=x$, then we can modify $T$ so that all the arcs that have been added enter $x$. Then $x$ is again a feed vertex of the new tournament $T'$. Since all the out-neighbors of $x$ in $T'$ are whole vertices (a whole vertex is avertex adjacent to all the other vertices in $G$)in $G$, then $N^{+}(x)=N^{+}_{T'}(x)$ and $N^{++}(x)=N^{++}_{T'}(x)$. However, $x$ has the SNP in $T'$, since it is a feed vertex, then $x$ has the SNP in $G$ as well. Finally, suppose that $v=y$ and $xy$ is a missing edge. Again we reorient the arc incident $y$ so that it enters $y$. Then $y$ is a feed vertex of the new tournament $T'$, and so $y$ has the SNP in $T'$. We have $N^{+}(y)=N^{+}_{T'}(y)$. Moreover, $N^{++}(y)=N^{++}_{T'}(y)$. Indeed, suppose that $(y,z), (z,a), (a, y)\in E(T')$. Then $z\neq x$. If $a\neq x$ or $z$ is a whole vertex, then $(y,z), (z,a)\in E$. Otherwise, $a=x$ and $z$ is not a whole vertex. So $zx$ is a missing edge and $(z,x)$ is an added arc to $G$. Then for the  missing edge $zx$ we have $R=\phi$. We have $(y,z)\in E$ and $(x,y)\notin E$ (because $xy$ is a missing edge of $G$). Then $x\in N^{++}(y)$, since otherwise, $y\in R=\phi$ which is a contradiction. Therefore, $y$ received no new vertex in its second out-neighborhood in $T'$. Thus $y$ has the SNP in $G$ as well.
  Similarly, if $v$ is a whole vertex, then we have $N^{++}(v)=N^{++}_T(v)$. Thus $v$ satisfies the SNP in $T$ and $G$ as well.


\begin{thebibliography}{}
\bibitem{m.o.}
F. Havet and S. Thomass\'{e}: Median Orders of Tournaments: A Tool for the Second Neighborhood
Problem and Sumner's Conjecture. J. Graph Theory 35, 244-256  (2000)

\bibitem{fidler}
D. Fidler and R. Yuster: Remarks on the second neighborhood problem. J. Graph Theory 55, 208-220 (2007)

%
%


\end{thebibliography}
\end{document}